\documentclass[12pt,draft,onecolumn]{IEEEtran}
\usepackage[final]{graphicx}
\usepackage{url}
\usepackage{color}
\usepackage{amsmath}
\usepackage{amssymb}
\usepackage{amsthm}
\usepackage{amsfonts}
\usepackage{bm}
\usepackage{psfrag}
\usepackage[mathscr]{eucal}
\usepackage{cancel}
\usepackage{times}
\usepackage{epsfig,color,pstricks}
\usepackage{pst-node,pst-blur}
\usepackage{amsmath,amssymb,amsthm,amsfonts,eucal,latexsym,mathcom}
\usepackage{cite}
\usepackage{multirow}
\usepackage{algorithm}
\usepackage{algorithmic}
\usepackage{graphicx}
\usepackage{makeidx}
\usepackage{subfig}
\usepackage{etoolbox}
\usepackage{enumitem}

\newtheorem{mydef}{Definition}

\newcommand{\bl}{\boldsymbol}
\newcommand{\bd}{\boldsymbol}

\newcommand{\pd}[2]{\frac{\partial #1}{\partial #2}}
\newtheorem{theorem}{Theorem}
\newtheorem{lemma}[theorem]{Lemma}

\newenvironment{remark}[1][Remark]{\begin{trivlist}
\item[\hskip \labelsep {\bfseries #1}]}{\end{trivlist}}

\title{A Geometric Analysis of Power System \\ Loadability Regions}
\author{{Yang~Weng,~\IEEEmembership{Member,~IEEE},}
{Ram Rajagopal,~\IEEEmembership{Member,~IEEE},}
{and Baosen Zhang,~\IEEEmembership{Member,~IEEE}}
\vspace{-11mm}}

\begin{document}

\maketitle
\begin{abstract}
Understanding the feasible power flow region is of central importance to power system analysis. In this paper, we propose a geometric view of the power system loadability problem. By using rectangular coordinates for complex voltages, we provide an integrated geometric understanding of active and reactive power flow equations on loadability boundaries. Based on such an understanding, we develop a linear programming framework to $1$) verify if an operating point is on the loadability boundary, $2$) compute the margin of an operating point to the loadability boundary, and $3$) calculate a loadability boundary point of any direction. The proposed method is computationally more efficient than existing methods since it does not require solving nonlinear optimization problems or calculating the eigenvalues of the power flow Jacobian. Standard IEEE test cases demonstrate the capability of the new method compared to the current state-of-the-art methods. 
\end{abstract}

\vspace{-4mm}
\section{Introduction}

The loadability of electrical networks considers whether operating points are feasible under physical constraints and its study has long been an integral part of power systems  planning and operation. In the planning phase, loadability analysis can be used to determine the need for shunt compensation, new transmission lines \cite{Shand24}, reserves \cite{Lajda81},  
and other system additions \cite{Johnson77}. In the operation phase \cite{Cutsem91}, a load flow solution can be used to find stability margins in preventing voltage collapses caused by large load variation and saddle-node bifurcation \cite{CutsemEtAl1998}. 
As more renewable resources are integrated into the aging electrical infrastructure and systems operate closer to their limits, characterizing the loadability of power systems is becoming increasingly important~\cite{Molzahn13}.

The most well-known example of a loadability limit is the power-voltage (PV) curve for a $2$-bus system, where the maximum loading occurs at the tip of the curve~\cite{McCalley12}. 
For larger systems, visualizing and computing the loadability boundary becomes more difficult. A typical approach is to increase the load on all buses by the same factor from a given base load until a power flow solution cannot be found~\cite{Zeng93}. 
However, as the load increases, the power flow Jacobian becomes ill-conditioned and standard first and second order algorithms become numerically unstable~\cite{crow2015}.

Overcoming the computational challenges at near the loadability limits has received considerable attention from the community. The work in  \cite{Ajjarapu92} develop a robust continuation power flow method for obtaining solutions on the P-V curve and \cite{GhiocelEtAl2014} propose an angle-reactive-power (AQ) bus to mitigate the numerical issues. A notable subset of these studies are the non-divergent power flow methods~\cite{Bijwe03,
Schaffer88}. In these methods, a voltage update increment is adjusted via computing a multiplier to avoid divergence while minimizing the norm of voltage error residuals. Other algebraic methods have been used to more directly study the properties of power flow solutions.  A popular approach is to examine the eigenvalues of the power flow Jacobian~\cite{Canizares93}, although the singularity of the Jacobian is necessary but not sufficient to conclude that a solution is on the boundary. A set of studies in~\cite{Cutsem95,Capitanescu05,Hamon13,Perninge11} research on how to use local information to understand whether a point is operating on the boundary or close to it.
Finally, polynomial homotopy continuation methods can be employed to find all the power flow solutions, although often at prohibitively high computational costs~\cite{Mehta14}.

These developments have proven to be quite useful in practice, but three challenges still remain. First, the geometric intuition of the loadability in $2$-bus PV curve needs to be extended to larger systems for loadability calculation. Second, most of the existing methods increase loads at all the buses by the same factor, thus limiting the exploration of the full region. Third, since nonlinear optimization problem is often employed, the computational requirements are nontrivial.

In this paper, we overcome these three challenges by $1$) extending the geometrical loadability intuition in the PV curve to more than $2$ buses and to include reactive power, and $2$) presenting a linear programming approach to test whether an operating point is on the loadability boundary, to find the distance between the point and the boundary, and to locate any loadability boundary point of interest. The starting point of our analysis is using rectangular coordinates to represent complex voltages and studying the power flow equations and Jacobian matrices~\cite{WeiEtAl1998,TorresEtAl1998,MolzahnEtAl2014}.

Formally, the loadability boundary of a system is the set of points on the Pareto-Front of active powers~\cite{zhang2013geometry}, where the load on one bus cannot be strictly increased without decreasing the load at some other buses while satisfying the physical constraints. Therefore, this Pareto-Front represents the limit of operating a system. We make the following contributions:
\begin{enumerate}
  \item From rectangular coordinates, we develop a geometric view for systems with more than $2$ buses that integrates real and reactive powers, motivating the analysis later on.
  \item We show that the eigenvalues of the power flow Jacobian are insufficient to describe system loadability boundary. Instead, we present a linear programming approach to test whether an operating point is on the boundary.
  \item Based on the linear programming approach, we characterize the \emph{loadability margins} of operating points. We also formulate a linear programming problem to characterize the power flow feasibility boundary points.
\end{enumerate}
We validate our approaches by simulations on different transmission grids such as the  $14$, $300$, and $13659$-bus networks and two distribution grids (the $8$ and $123$-bus networks~\cite{Kersting01}).

The rest of the paper is organized as follows: Section~\ref{sec:model} motivates the rectangular coordinate-based analysis and provides an integrated geometric view of active/reactive power flow equations. Section~\ref{sec:jacobian} shows the linear Jacobian matrix and its application for security boundary point verification. Section~\ref{sec:marginal} quantifies margins of points that are not on the boundary. Section~\ref{sec:pareto} shows how to search for all of the boundary points. Section~\ref{sec:num} evaluates the performance of the new method and Section~\ref{sec:con} concludes the paper.

\section{Rectangular Coordinates and Geometry of Power Flow}\label{sec:model}
\subsection{Visualization of Complex Power Flow}
The power flow equations in polar coordinates are \cite{crow2015}:
\begin{subequations}
\label{equ:pfqf_polar}
\begin{align}
p_d = &\sum_{k=1}^{n}|v_d||v_k|\left(g_{dk}\cos\theta_{dk} + b_{dk}\sin\theta_{dk}\right),\label{equ:pf_polar}\\
q_d = &\sum_{k=1}^{n}|v_d||v_k|\left(g_{dk}\sin\theta_{dk} - b_{dk}\cos\theta_{dk}\right),\label{equ:qf_polar}
\end{align}
\end{subequations}
where $n$ is the number of buses in the network; $p_d$ and $q_d$ are the active and reactive power injections at bus $d$; $v_d$ is the complex phasor at bus $d$ and $|v_d|$ is its voltage magnitude; $\theta_{dk}=\theta_k-\theta_d$ is the phase angle difference between bus $k$ and bus $d$; $g_{dk}$ and $b_{dk}$ are the electrical conductance and susceptance between bus $d$ and bus $k$. Together, $y_{dk} = g_{dk}+j\cdot b_{dk}$ forms the admittance, where $j$ is the imaginary unit.

These equations have been the central objects of interest in power system analysis for decades. Because of the nonlinear interaction of sinusoidal and polynomial functions, the power flow-based loadability analysis remains challenging. For some simple cases, such as a $2$-bus system, PV curve visualizes the system loadability. 
However, this geometric picture is hard to extend while keeping all system information for joint loadability analysis. 
For example, the PV curve has been extended to large systems by multiplying each bus by a loading factor, then visualizing the impact on voltage stability as this loading factor changes~\cite{CutsemEtAl1998,Perninge11}. But this approach hides the local behavior of each bus and picking a good starting load is not always easy.

\subsection{Rectangular Coordinate-based Power Flow}
One difficulty in~\eqref{equ:pfqf_polar} for loadability analysis lies in its diverse functional types, e.g., sinusoidal and polynomial. To reduce the functional types for easier loadability analysis, we adopt the rectangular coordinates for complex voltages.
Let $v_{d,r} \triangleq Re(v_d)$ and $v_{d,i} \triangleq Im(v_d)$ be the real and imaginary parts of the complex voltage at bus $d$, respectively. Then, the power flow equations in \eqref{equ:pfqf_polar} become
\begin{subequations}\label{equ:pfqf_rec}
\begin{align}
p_{d} =t_{d,1}\cdot v_{d,r}^2 + t_{d,2}\cdot v_{d,r}+ t_{d,1}\cdot v_{d,i}^2 + t_{d,3}\cdot v_{d,i}\label{equ:pf_rec},\\
q_{d} =t_{d,4}\cdot v_{d,r}^2 - t_{d,3}\cdot v_{d,r}+ t_{d,4}\cdot v_{d,i}^2 + t_{d,2}\cdot v_{d,i}\label{equ:qf_rec},
\end{align}
\end{subequations}
in rectangular coordinates, where
\begin{subequations}\label{equ:t_def}
\begin{align}
&t_{d,1} = -\sum_{k\in \mathcal{N}(d)}g_{kd},~ t_{d,2} = \sum_{k\in \mathcal{N}(d)}(v_{k,r}g_{kd}- v_{k,i}b_{kd}),\\
&t_{d,3} = \sum_{k\in \mathcal{N}(d)}(v_{k,r}b_{kd}+ v_{k,i}g_{kd}),~t_{d,4} = \sum_{k\in \mathcal{N}(d)}b_{kd},
\end{align}
\end{subequations}
where $\mathcal{N}(d)$ is the neighbors of bus $d$. The detailed derivations from \eqref{equ:pfqf_polar} to \eqref{equ:pfqf_rec} are given in the Appendix \ref{app:der}. 

One benefit of using  \eqref{equ:pfqf_rec} comes from its capability of visualizing active and reactive power flow equations in the same space. Such an understanding will directly give intuitive meaning of Pareto Front for the loadability boundary considered in the next section. Specifically, for fixed constants $t_{d,1},t_{d,2},t_{d,3},t_{d,4}$, \eqref{equ:pf_rec} and \eqref{equ:qf_rec} describe two circles in the $v_{d,r}$ and $v_{d,i}$ space. Note that the circle concept is different than the power circle concept in the past \cite{Chard54,Goodrich52}. The next lemma characterizes the centers and radii of the active and reactive power flow circles.
\begin{lemma} \label{lem:circle}
The centers and radii: The coordinates of circle center $E$ for the active power flow are $\left(-\frac{t_{d,2}}{2t_{d,1}}, -\frac{t_{d,3}}{2t_{d,1}}\right)$ for bus $d$. Its radius decreases when $p_{d}$ increases. The coordinates of circle center $D$ for the reactive power flow are $\left(\frac{t_{d,3}}{2t_{d,4}}, -\frac{t_{d,2}}{2t_{d,4}}\right)$ for bus $d$. Its radius decreases when $q_{d}$ increases.
\end{lemma}
\begin{proof}
See Appendix \ref{lem:circle_app}.
\end{proof}

\begin{figure*}[!ht]
  \centering
    \includegraphics[width=0.9\linewidth]{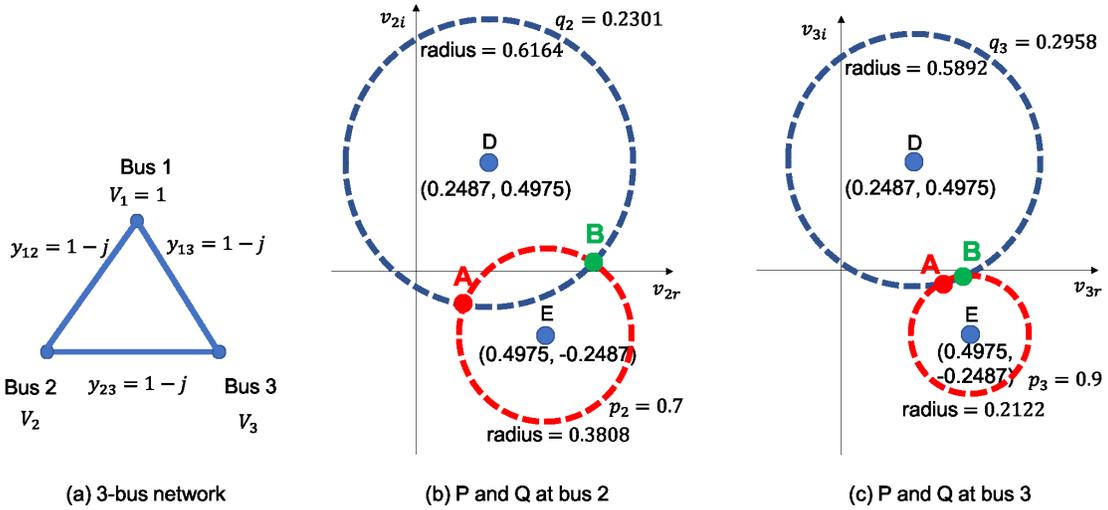}
  \caption{A three bus network and the circles formed at buses $2$ and $3$. The closeness of the intersections at bus $3$ indicates that the system is operating close to the boundary of the feasible region.}
  \label{fig:tri_all}
  \end{figure*}

For example, consider the $3$-bus network in Fig.~\ref{fig:tri_all}, where bus 1 is the slack bus, and bus $2$ and $3$ are PQ load buses. Let the admittance be $1-0.5j$ for all the lines, $p_2=0.7$, $p_3=0.9$, power factor$=0.95$. Fig. \ref{fig:tri_all}(b) and Fig. \ref{fig:tri_all}(c) show the circles formed by \eqref{equ:pfqf_rec} for bus $2$ and bus $3$, respectively. In this case, the intersection points between the two circles in Fig~\ref{fig:tri_all}(b) are far apart, whereas the two points on the intersection of the two circles in Fig.~\ref{fig:tri_all}(c) are close together. This suggests that the system is operating close to its limit and small changes may lead to insolvability of the power flow equations at bus $3$. Notably, Fig~\ref{fig:tri_all}(b) and Fig~\ref{fig:tri_all}(c) is a generalization of the well-studied concept of PV curve for power systems. 

\section{The Boundary of Power Flow Region: Beyond the Jacobian}\label{sec:jacobian}
When the system is approaching the loadability boundary, e.g., nose point of a $2$-bus system, the point A and point B will come closer. Therefore, we  need to algebraically characterize the operating points, 
especially for those on or close to the boundary of the feasible power flow region. Commonly, these types of analysis are done through the power flow Jacobian, and in particular, the singularity of the Jacobian matrix has long been used to characterize the solvability and stability of power flow solutions~\cite{KundurEtAl1994a}. In this section, we show that the Jacobian is \emph{not always sufficient to identify the boundary of the power flow region}, and we propose a different method of identifying whether a solution is on the boundary by using a linear program in the rectangular coordinates.

\subsection{The Limitation of Singularity Analysis of Jacobian Matrix}
The power flow Jacobian matrix, $\bl J$, is normally defined by the first order partial derivatives of active and reactive powers with respect to the state variables. In our analysis, these partial derivatives are taken with respect to the real and imaginary parts of bus voltages:
\begin{equation*}
	\bl J= \begin{bmatrix} \pd{\bl p}{\bl v_{r}} & \pd{\bl p}{\bl v_i} \\
	\pd{\bl q}{\bl v_r} & \pd{\bl q}{\bl v_i} \end{bmatrix},
\end{equation*}
where the elements are given by
\begin{subequations} \label{eqn:partial_rec}
\begin{equation}
\frac{\partial p_d}{\partial v_{k,r}} = \left\{
\begin{array}{ll}
2t_{d,1}v_{d,r}+t_{d,2},~~& \text{if } d = k,\\
g_{kd}v_{d,r}+b_{kd}v_{d,i},~~& \text{if } d \neq k.
\end{array} \right.
\end{equation}
\begin{equation}
\frac{\partial p_d}{\partial v_{k,i}} = \left\{
\begin{array}{ll}
2t_{d,1}v_{d,i}+t_{d,3}, & \text{if } d = k,\\
-b_{kd}v_{d,r}+g_{kd}v_{d,i}, & \text{if } d \neq k.
\end{array} \right.
\end{equation}
\begin{equation}
\frac{\partial q_d}{\partial v_{k,r}} = \left\{
\begin{array}{ll}
2t_{d,4} v_{d,r}-t_{d,3}, & \text{if } d = k,\\
-b_{kd}v_{d,r}+g_{kd}v_{d,i}, & \text{if } d \neq k.
\end{array} \right.
\end{equation}
\begin{equation}
\frac{\partial q_d}{\partial v_{k,i}} = \left\{
\begin{array}{ll}
2t_{d,4} v_{d,i}+t_{d,2}, & \text{if } d = k,\\
-g_{kd}v_{d,r}-b_{kd}v_{d,i}, & \text{if } d \neq k.
\end{array} \right.
\end{equation}
\end{subequations}
For the partial derivatives above, $t_{d,1}$ and $t_{d,4}$ are both constant given network parameters, and $t_{d,2}$ and $t_{d,3}$ are \emph{linear} in the variables. Therefore, each of the partial derivatives in \eqref{eqn:partial_rec} is linear in the state variables $v_{d,r}$ and $v_{d,i}$.

The Jacobian is normally used via the inverse function theorem, which states that the power flow equations stabilize around an operating voltage if the Jacobian is non-singular. This condition is necessary since every stable point must have a nonsingular Jacobian. However, the singularity of the Jacobian is insufficient~\cite{Bijwe03}, especially for the loadability boundary. As the next example will show, a singular Jacobian does not imply that the operating voltage is at the boundary of the power flow feasibility region.

\begin{figure}[ht!]
	\centering
	\includegraphics[scale=0.4]{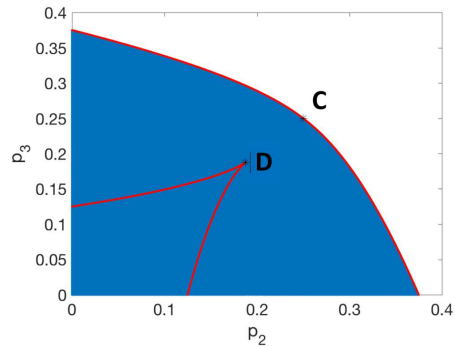}
	\caption{Feasible power flow region for the $3$-bus network. The red points denote operating points when the Jacobian is singular. They separate into two parts: the boundary of the region (e.g., point F) and points in the strict interior (e.g., point H).}
	\label{fig:jacobian}
\end{figure}

Again, consider the $3$-bus network in Fig.~\ref{fig:tri_all}. For simplicity, we assume the lines are purely resistive (all line admittances are $1$ per unit) and only consider active powers. Let bus $1$ be the slack bus and buses $2$ and $3$ be load buses consuming positive amount of active powers.
In this case, the Jacobian becomes:
\begin{equation}\label{eqn:Jacobian_tri}
	\bl J=\begin{bmatrix} 1-4v_2+v_3 & v_2 \\ v_3 & 1-4v_3+v_2 \end{bmatrix},
\end{equation}
where $v_2$ and $v_3$ are the voltages at bus $2$ and bus $3$, respectively. Fig. \ref{fig:jacobian} shows the feasible power flow region of power consumptions at buses $2$ and $3$. The red lines show the points where the Jacobian is singular.

Here, we focus on two particular points in Fig.~\ref{fig:jacobian}, points $F$ and $H$. 
At these points, $v_1=v_2=v$ due to the symmetry of the network. Then, finding the determinant of $J$ and equating it to $0$, we obtain $v=0.25$ (point $H$) or $v=0.5$ (point $F$). We emphasize that these two points are qualitatively different. Point $F$ is on the boundary of the feasible region, and therefore is a loadability point. However, point $H$ is well within the strict interior of the ``feasible region", therefore it is not on the loadability boundary. 
Points like $H$ are sometimes called cusp bifurcation points in stability analysis \cite{Vournas99,Harlim07}. Therefore, if we are interested in finding whether a point is on the boundary or characterizing its loadability margin, just looking at the determinant of the Jacobian is insufficient.

\begin{figure}[ht]
    \centering
    \includegraphics[scale=0.45]{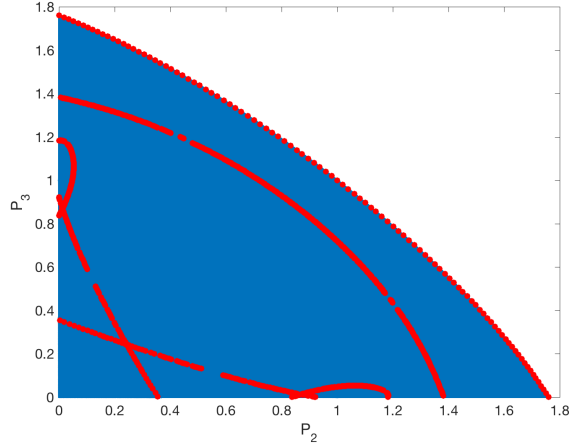}
    \caption{The active power injections where the Jacobian is singular for a $3$-bus fully connected network, with reactance $j$ per unit on each of the lines. There are two distinct types of points (marked in red), the boundary and the interior, that both have a singular Jacobian matrix.}
    \label{fig:Jacobian2}
  \end{figure}

In addition to the purely resistive network to generate Fig.~\ref{fig:jacobian}, we also plot the feasibility region and the points where the Jacobian is singular for a network with complex impedances as in~Fig.~\ref{fig:Jacobian2}. Other cases have results similar to Fig.~\ref{fig:jacobian} and Fig.~\ref{fig:Jacobian2}. 

\subsection{Verify Loadability Points via Pareto-Front Method}

To isolate just the points on the boundary, we need to look deeper into the power flow equations than just the determinant of the Jacobian. Here, we focus on a network where the buses are loads. Geometrically, a point is on the \emph{loadability boundary} (or simply boundary) if there does not exist another point that can consume more power:
\begin{mydef}\label{def:boundary}
	Let $\bl v$ be the complex voltages and $\bl p=(p_1,\cdots,p_n)$ be the corresponding bus active powers. We say that the operating point is on the loadability boundary if there does not exist another operating point $\hat{\bl p}=(\hat{p}_1,\cdots,\hat{p}_n)$ such that $\hat{p}_k \geq p_k$ for all bus index $k$ (nonnegative change in load) and $\hat{p}_d > p_d$ for at least one $d\neq k$ (positive change in at least one load).
\end{mydef}
This definition coincides with the definition of the Pareto-Front since we are modeling each bus (except the slack) as a load bus. It can be easily extended to a network where some buses are generators by changing the direction of inequalities in the definition above. Instead of looking at the determinant of the Jacobian, the next theorem gives a linear programming condition for the points on the boundary:
\begin{theorem}~\label{thm:boundary}
	Checking whether an operating point is on the boundary of the feasible power flow region is equivalent to solving a linear programming problem, e.g., \eqref{opt:verify}.
\end{theorem}


Let 
\begin{align}
\bl h_d=[\pd{p_d}{v_{1,r}}, \pd{p_d}{v_{1,i}}, \pd{p_d}{v_{2,r}}, \cdots , \pd{p_d}{v_{n,i}}]^T \in R^{2n}
\end{align}
be the gradient of $p_d$ with respect to all the state variables. Therefore, $\bl h_d$ is the transpose of the $d^{th}$ row of the Jacobian matrix. Let $\bl z \in R^{2n}$ be a direction, towards which we move the real and imaginary part of the voltages. Then, by Definition~\ref{def:boundary}, a point is on the boundary if there does not exist a direction to move where the consumption of one bus is increased without decreasing the consumption at other buses.

Therefore, we can check if there is a direction $\bl y$ that makes the following problem feasible. Suppose that $\bl h_1,\cdots,\bl h_n$ are given.
\begin{subequations} \label{opt:feasible}
\begin{align}
	& \bl y^T \bl h_d \geq 0, \mbox{ for all } d=1,\cdots,n ,\label{eqn:verify_all} \\
	& \sum_{d=1}^n \bl y^T \bl h_d =1. \label{eqn:verify_1}
\end{align}
\end{subequations}
The constraint \eqref{eqn:verify_all} specifies that moving in the direction $\bl y$ cannot decrease any of the active powers. The constraint \eqref{eqn:verify_1} is equivalent to stating that at least one bus' active power must strictly increase. This comes from the fact that $\bl y$ is not a constraint. Therefore, as long as $\bl y^T\bl h_d >0$ for some $d$, the sum $\sum_{d=1}^n \bl y^T \bl h_d$ can be scaled to be $1$. If the problem (\ref{opt:feasible}) is feasible, the corresponding power pair is not on the boundary. If the problem (\ref{opt:feasible}) is infeasible, this means that the point is on the Pareto-Front. So, it is on the loadability boundary.

By adding a constant objective, we can encode this condition in a linear programming (LP) feasibility problem. This is because, in a constraint optimization problem, a solver usually tries to firstly find a feasible set by using the constraints. Then, it will use searching methods, e.g., gradient descent method, for the objective in the feasible region. Therefore, by converting the feasibility problem (\ref{opt:feasible}) into an optimization form, we can use the state-of-the-art solver in convex optimization tool set, which is quite efficient.

Then, we solve the following:
\begin{subequations} \label{opt:verify}
\begin{align}
	\min_{\bl y} \;\; & 1 \\
	\mbox{s.t. } & \bl y^T \bl h_d \geq 0, \mbox{ for all } d=1,\cdots,n ,\label{eqn:verify_all_2} \\
	& \sum_{d=1}^n \bl y^T \bl h_d =1. \label{eqn:verify_1_2}
\end{align}
\end{subequations}
In this optimization problem, the objective is irrelevant since we are only interested in whether the problem is feasible. Finally, an operating point is on the boundary if and only if the problem in \eqref{opt:verify} is infeasible.

A system operating on the boundary limit will lose stability before our conditions are checked. So, we provide an alarm when a system is approaching this boundary for practical interest. In the following, we change the optimization (\ref{opt:verify}) slightly to provide an alarm by setting up an $\epsilon$ value for earlier alarming.   
\begin{subequations} \label{opt:verify2}
\begin{align}
	\min_{\bl y} \;\; & 1 \\
	\mbox{s.t. } & \bl y^T \bl h_d \geq \epsilon, \mbox{ for all } d=1,\cdots,n ,\label{eqn:verify_all_3} \\
	& \sum_{d=1}^n \bl y^T \bl h_d =1. \label{eqn:verify_1_3}
\end{align}
\end{subequations}

\begin{figure*}[ht]
  \centering
  \subfloat[Margin according to Thevenin method.]{ \label{fig:margin_th}
    \includegraphics[width=0.45\textwidth]{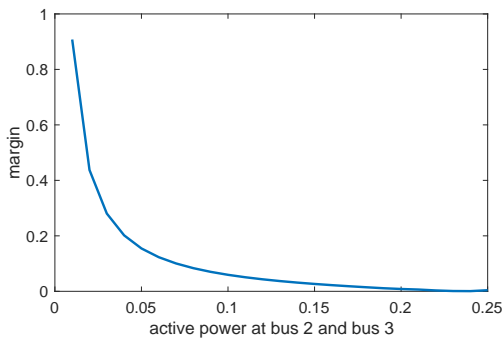}}
\subfloat[Margin using the proposed method.]{ \label{fig:margin}
  \includegraphics[width=0.45\textwidth]{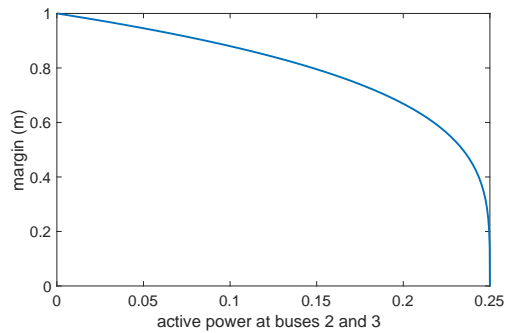}}
  \caption{Comparison of the margin computed by (a) the Thevenin equivalent method and (b) our proposed method.}
  \label{fig:3-bus}
\end{figure*}

\begin{remark}
	In the system operating on or near the loadability boundary, computation speed is important, otherwise the system may lose stability before we can compute anything. The algorithm we propose involves solving a simple linear program that is on the order of the system size, e.g., bus numbers, making the computational speed very short even for large systems. In contrast, some of other nonlinear calculators require iterative methods that solve successive nonlinear problems, resulting in a much slower process.
\end{remark}

For the example in Fig.~\ref{fig:jacobian}, point $A$ is given by $v_2=v_3=0.5$ and point $B$ is given by $v_2=v_3=0.25$. Using \eqref{eqn:partial_rec}, we have $\bl h_2=[-0.5, 0.5]^T$ and $\bl h_3=[0.5, -0.5]^T$ at point $A$, and $\bl h_2=\bl h_3=[0.25, 0.25]^T$ at point $B$. It is easy to check that a $\bl y$ feasible for \eqref{opt:verify} exists at point $B$. However, it is impossible to find such a $\bl y$ at point $A$. Therefore, we can conclude that $A$ is on the boundary whereas $B$ is not, even though both have singular Jacobians.
For how to add more constraint, please refer to Appendix \ref{sec:add}.

\section{Loadability Margins}
\label{sec:marginal}
In addition to asking if a point is on the boundary, we are sometimes more interested in \emph{how close} a point is to the boundary. For this purpose, we will apply a simple modification to \eqref{opt:verify} for measuring the distance, or the margin.

\subsection{Measure the Loadability Margin}
In the optimization problem \eqref{opt:verify}, we check \emph{if it is possible} to move the operating point in a direction such that the active powers can be increased. For a point close to the boundary, we are interested in \emph{how much} a point can be moved before reaching the boundary. Therefore, we use the optimal objective value in \eqref{eqn:margin} to measure the stability \emph{margin} of an operating point. Again, let $\bl h_d^T$ be the $d^{th}$ row of the Jacobian matrix at an operating point of interest. Then, we solve

\begin{subequations}\label{eqn:margin}
  \begin{align}
    m=\max_{\bl y} & \sum_{d=1}^n \bl y^T \bl h_d \\
    \mbox{s.t. } & \bl y^T \bl h_d \geq 0, \mbox{ for all } d=1,\cdots,n, \\
    & ||\bl y||_2 = 1. \label{eqn:norm}
  \end{align}
\end{subequations}
In this problem, we look for a unit vector $\bl y$ such that the sum of the active powers can be increased by the maximum amount. The value of the optimization problem is denoted by $m$, which we think as the margin or the distance to the boundary. Note that the constraint in \eqref{eqn:norm} may seem non-convex, but for any point that is not on the boundary, \eqref{eqn:norm} can be relaxed to $||\bl y||_2 \leq 1$ without changing the objective value. Therefore, \eqref{eqn:margin} can be easily solved by standard solvers.


\subsection{A Comparison to Thevenin-Equivalent Margin}
Here, we compare the solution of \eqref{eqn:margin} with the solution of widely adopted Thevenin-equivalent method for margin calculation \cite{Vu99}. 
Specifically, \cite{Vu99} describes a technique to find the margin and condition for maximum loadability. The bus of interest is considered as the load bus and the rest of the system is replaced with a Thevenin impedance and Thevenin voltage. Originally, there are two voltage solutions for a given power transfer. As the power transfer reaches its maximum value, there is only one voltage solution and this point is known as bifurcation point. Mathematically, Kirchhoff's current law leads to $(p+jq)\cdot v^{*}_{Thev} = v\cdot ({v} - {e})^{*}$, where $e$ is the voltage at the aggregated infinite bus and $(\cdot)^{*}$ is the complex conjugate operator. At the system bifurcation point, $v = (e - v)^{*}$. Therefore, $z_{App}\cdot i
   = (z_{Thev}\cdot i)^{*}$. This leads to $|{z}_{App}| = |{z}_{Thev}|$.

Hence, by tracking how close the Thevenin impedance is to the load impedance, we can know the margin for the maximum loading condition. When ${z}_{Thev} = {z}_{app}$, there is only one voltage solution and hence the maximum power transfer capacity is reached. For ease of illustration, we adopt the $3$-bus system with bus $1$ being the generator (and slack), buses $2$ and $3$ being loads. We set the active power at buses $2$ and $3$ to be equal and increase it until the system becomes unstable. The resulting margin is shown in Fig.~\ref{fig:3-bus}.

From Fig.~\ref{fig:3-bus}, we see that both methods show that the system has a margin of $0$ when the load is at $0.25$ p.u. However, the Thevenin equivalent method is much more conservative than ours. For example, when the load is half of the maximum load ($0.125$ p.u.), the Thevenin equivalent method has a relatively small margin, which may lead an operator to conclude that the load cannot be increased much more and operate conservatively. This would result in inefficiencies in operations, especially in an aging grid that is facing more complex loading environments~\cite{VaimanEtAl2010}.  In contrast, our method provides a much larger margin when the load is far away from the maximum, and the margin decreases rapidly once the load approaches the maximum. This allows operators to better gauge the state of the system, leading to more efficient and reliable operations.

\section{Locating All Lodability Boundary Points via Pareto-Front}
\label{sec:pareto}

In the last two sections, we have explored how to determine whether a point is on the Pareto-Front and its ``distance'' to the front. In this section, we ask the question of whether we can determine the Pareto-Front itself. To answer this question, we observe that by definition, for any point $\bl p$ on the Pareto-Front, there exists a vector $\bl z$ such that $\bl z^T \bl p$ is the maximum among all possible active power vectors. Conversely, by varying $\bl z$ and maximizing over $\bl p$, we can find the Pareto-Front. Here, $\bl z$ physically means the direction that the powers at different buses are growing. Depending on the direction or the ratio of loads on different buses, one can find a boundary at certain loading ratio conditions.

Of course, for a large system, exhaustively varying $\bl z$ is impractical. However, in many cases, there are a few $\bl z$'s of special interest. For example, if $\bl z$ is the all-ones vector $\bl 1$, we are looking for the maximum sum power that the network can support. In other settings, there are a few classes of loads, and $\bl z$ has only a few distinct values.

To find an active power vector $\bl p$ such that $\bl z^T \bl p$ is the maximum, we again look at the partial derivatives of active power with respect to the voltages.
On the Pareto-Front and given a $\bl z$ tangent to it, the gradients of active power with respect to voltages are \emph{orthogonal} to $\bl z$. Let $\bl h_d^T$ be the $d^{th}$ row of the Jacobian, which we decompose into two parts: $\bl h_{d,r}^T$ is the first $n$ components corresponding to the $d^{th}$ row of $\frac{\partial \bl p}{\partial \bl v_r}$ and $\bl h_{d,i}^T$ is the last $n$ components corresponding to the $d^{th}$ row of $\frac{\partial \bl p}{\partial \bl v_i}$. We then look for operating points such that $\bl z^T \bl h_{d,r}=0$ and $\bl z^T \bl h_{d,i}=0$ for all $d=1,\cdots,n$.

Rectangular coordinates make this problem much easier to solve. As shown in \eqref{eqn:partial_rec}, the elements of the Jacobian are \emph{linear} in the real and imaginary parts of the bus voltages. Therefore, for a given $\bl z$, the system of equations $\bl z^T \bl h_{d,r}=0$ and $\bl z^T \bl h_{d,i}=0$ for all $d=1,\cdots,n$ becomes \emph{a system of linear equations}:
\begin{subequations}\label{equ:bus_d}
  \begin{align}
    \left(2t_{d,1}v_{d,r}+t_{d,2}\right)z_d&+\sum_{k\neq d}\left(g_{kd}v_{d,r}+b_{kd}v_{d,i}\right)z_k=0,\\
    \left(2t_{d,1}v_{d,i}+t_{d,3}\right)z_d&+\sum_{k\neq d}\left(g_{kd}v_{d,i}-b_{kd}v_{d,r}\right)z_k=0,\\
    d&=1,\cdots,n. \nonumber
  \end{align}
\end{subequations}
Recall that $t_{d,1}=-\sum_{k\in \mathcal{N}(d)}g_{kd}$ and every equation is linear in $v_{d,r}$ and $v_{d,i}$ once $\bl z$ and the network are given.

\begin{figure}[ht]
    \begin{center}    \includegraphics[width=0.5\textwidth]{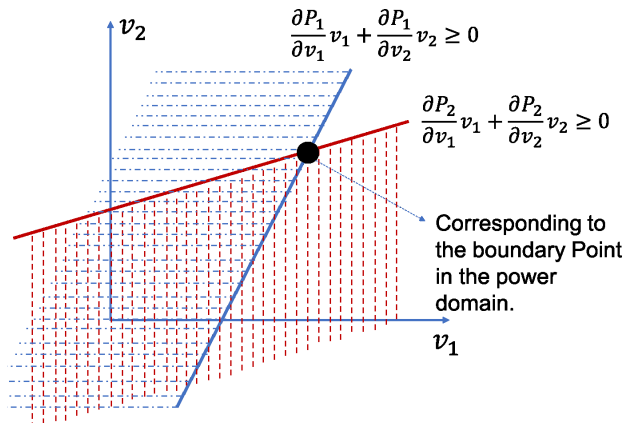}
    \caption{A point on the Pareto-Front for the three bus network in Fig.~\ref{fig:tri_all} obtained from solving \eqref{equ:bus_d} with $\bl{z}=[1\;1]^T$. This point maximizes the sum power $P_1+P_2$. }
    \label{fig:boundarypoint}
    \end{center}
\end{figure}

Fig. \ref{fig:boundarypoint} shows the solution of \eqref{equ:bus_d} for the $3$-bus network in Fig.~\ref{fig:tri_all} with all line admittances being $1$'s and $\bl{z}=[1, 1]^T$. In this case, the network is purely resistive and its Jacobian is given in \eqref{eqn:Jacobian_tri}. Solving \eqref{equ:bus_d} gives $v_2=0.5$ and $v_3=0.5$, which corresponds to the point $A$ ($p_2 = 0.25$, $p_3 = 0.25$) on the boundary in Fig.~\ref{fig:jacobian}.

\begin{remark}
So far, we considered two bus types in our modeling: the reference bus and the load bus. We did not consider the generator bus type because the solar generator in the distribution grid can be modeled as a PQ bus \cite{Messenger10}.
\end{remark}

\section{Simulations}\label{sec:num}

\begin{table*}[ht]
\centering
\caption {Feasibility Boundary}
\label{tab:boundary}
\begin{tabular}{c*{14}{c}}
\hline
\hline
Bus No. & 3 & 4 & 5 & 6 & 9Q & 9target & 14 & 24 & 30 & 30pwl & 30Q & 39 & 57 & 89\\
\hline
On Boundary? & No & No & No & No & No & No & No & No & No & No & No & No & No & No\\
\hline
Time(s) & 1.5 & 1.6 & 1.6 & 1.6 & 1.6 & 1.4 & 1.7 & 1.9 & 1.6 & 2.0 & 1.8 & 1.7 & 1.6 & 1.7\\
\hline
\hline
Bus No. & 118 & 145 & 300 & 1354 & 1888 & 1951 & 2383 & 2736 & 2737 & 2746wop & 2746wp & 2848 & 2868 & 2869 \\
\hline
On Boundary? & No & No & No & No & No & No & No & No & No & No & No & No & No & No\\
\hline
Time(s) & 1.6 & 1.6 & 1.8 & 3.2 & 5.4 & 5.0 & 6.4 & 9.7 & 9.0 & 8.7 & 8.9 & 11.8 & 11 & 12.7 \\
\hline
\hline
Bus No. & 3012 & 3120 & 6468 & 6470 & 6495 & 6515 & 9241 & 13659 & 8 (Dist.) & 123 (Dist.)\\
\hline
On Boundary? & No & No & No & No & No & No & No & No & No & No\\
\hline
Time(s) & 11.7 & 12.4 & 45.5 & 50.6 & 42.4 & 44.8 & 102.1 & 281.9 & 1.4 & 1.5\\
\hline
\hline
\end{tabular}
\end{table*}

\begin{table*}[ht]
\centering
\caption {Loadability Margins via the Proposed Method}
\label{tab:margin}
\begin{tabular}{c*{14}{c}}
\hline
\hline
Bus No. & 3 & 4 & 5 & 6 & 9Q & 9target & 14 & 24 & 30 & 30pwl & 30Q & 39 & 57 & 89\\
\hline
Margin & 13.5 & 27 & 156.5 & 6.4 & 7.4  & 14.0 & 7.7 & 38.7 & 6 & 6 & 6 & 20  & 20 & 30\\
\hline
Time(s) & 3.3 & 3.3 & 3.4 & 3.4 & 3.3 & 3.3 & 3.4 & 3.3 & 3.3 & 3.6 & 3.5 & 3.4 & 3.3 & 3.4\\
\hline
\hline
Bus No. & 118 & 145 & 300 & 1354 & 1888 & 1951 & 2383 & 2736 & 2737 & 2746wop & 2746wp & 2848 & 2868 & 2869 \\
\hline
Margin & 8.6 & 1606 & 69.9 & 32.1 & 360.2 & 345.5 & 21.0 & 14.6 & 14.1 & 10.9 & 11.0 & 679.9  & 672.6 & 35.8\\
\hline
Time(s) & 3.5 & 3.5 & 3.7 & 5.1 & 7.5 & 7.4 & 10.7 & 10.5 & 10.7 & 10.4 & 10.4 & 11.3 & 11.6 & 12.7\\
\hline
\hline
Bus No.  & 3012 & 3120 & 6468 & 6470 & 6495 & 6515 & 9241 & 13659 & 8 (Dist.) & 123 (Dist.)\\
\hline
Margin & 16.4 & 15.7 & 1741.5 & 1769.3 & 1754.0  & 1756.6 & 23.1& 17.1 & 79 & 424\\
\hline
Time(s) & 11.7 & 12.4 & 45.5 & 50.6 & 42.4 & 44.8 & 102.1 & 281.9 & 1.5 & 2.0\\
\hline
\hline
\end{tabular}
\end{table*}

In this section, we perform extensive simulations on both transmission grids and distribution grids. The transmission grid cases include the standard IEEE transmission benchmarks ($4$, $9$, $14$, $30$, $39$, $57$, $118$, $300$-bus networks) and the MATPOWER test cases ($3$, $5$, $6$, $24$, $89$, $145$, $1354$, $1888$, $1951$, $2383$, $2736$, $2737$, $2746$, $2848$, $2868$, $2869$, $3012$, $3120$, $6468$, $6470$, $6495$, $6515$, $9241$, $13659$-bus networks)~\cite{Zimmerman09,Zimmerman10}. The distribution grids include standard IEEE $8$ and $123$-bus networks. The goal is to illustrate the three applications: $1$) verifying if a point is on the loadability boundary, $2$) measuring the loadability margins if the points are not on the boundary, and $3$) locating boundary points. 
As our proposed methods are based on linear programming and convex optimization, the computation time scales quite well with the size of the networks. The results on most of the test cases are similar to each other and we provide several representative examples in the followings.

\subsection{Boundary and Loadability Margins}
First of all, we use \eqref{eqn:partial_rec} to calculate the partial derivatives and then use \eqref{opt:verify} to test whether the operating points contained in the test cases are on the loadability boundary with Table \ref{tab:boundary}. Not surprisingly, none of the points are on the boundary as shown in Table~\ref{tab:boundary}. From the computation times in Table~\ref{tab:boundary}, we can see that for systems with thousands of buses, the condition in \eqref{opt:verify} can be checked in around $10$ seconds (using the CVX package in Matlab~\cite{cvx,gb08}) with an $i5$ laptop and $8$GB memory.


Next, we check the loadability margin of the operating points, recorded in Table~\ref{tab:margin}. Note that some systems are actually operating with fairly small margins, e.g., the $2746$-bus, $9241$-bus, and $13659$-bus systems. This means that they may not be robust under perturbations. The computational time again scales quite well with the size of the network. For comparison, we also list partial result of the Thevenin-Equivalent margin in Table \ref{tab:margin_TE}. The margin number is much smaller than the margin calculated from the proposed method in Table \ref{tab:margin}. So, our method provides a much larger margin when the load is far away from the maximum. This allows operators to better gauge system states for reliable operations.

\begin{table}[ht]
\centering
\caption {Thevenin-Equivalent Margin for Comparison}
\label{tab:margin_TE}
\begin{tabular}{c*{14}{c}}
\hline
\hline
Bus No. & 57 & 89 & 118 & 145 & 300 & 1354 & 1888 & 1951\\
\hline
Margin & 5.4 & 3.5 & 3.4 & 1.1 & 1  & 3 & 3 & 1.1\\
\hline
Time(s) & 2.5 & 2.1 & 2.6 & 2.2 & 2.6 & 31.12 & 68.9 & 8.3 \\
\hline
\hline
\end{tabular}
\end{table}

To show the impact of reactive power limits on the proposed method, we conduct a case study by adding a reactive power limit $-50<q_{69}<50$ at bus $69$ of the IEEE $118$-bus test case. 
Then, we compute the margin with such a reactive power limit via \eqref{opt:verify_q}. The margin is $6.1$. We also obtain the margin without such a reactive power limit, which is $8.6$. Therefore, the reactive power limit is binding at the operating point, which reduces the margin from $8.6$ to $6.1$.

\subsection{Going Towards the Boundary}

\begin{figure}[!h]
    \begin{center}
    \includegraphics[width=0.4\textwidth]{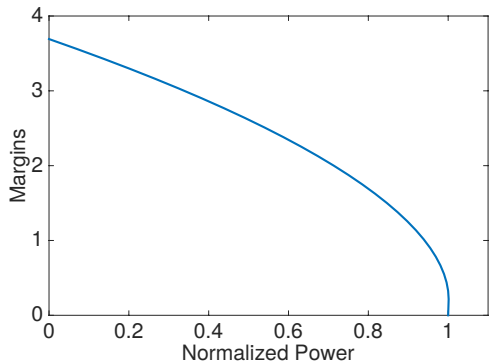}
    \caption{Margins vs. increasing active power for the IEEE $14$-bus system. Here, we successively increase the load on the buses (proportionally) to observe the change with respect to the margin.}
    \label{fig:towards}
    \end{center}
\end{figure}

Here, we move the operating points of a network in a direction until it hits the Pareto-Front. We first use the operating point calculated by running power flow from Matpower. Then, we use the obtained voltages, namely the operating point, as a direction in the state space to find the boundary point. Then, we change the voltage step-by-step from the Matpower-based voltage state to the boundary point that we obtain. The x-axis of Fig. \ref{fig:towards} shows the normalized active power as we successively increase the load on the IEEE $14$-bus system, and the y-axis plots the change in the margin.  As expected, the margin decreases when the point is moving towards the boundary point. When it is on the boundary point, the margin becomes zero.

\subsection{Locating Boundary Points}
\begin{figure}[!h]
    \begin{center}
    \includegraphics[width=0.4\textwidth]{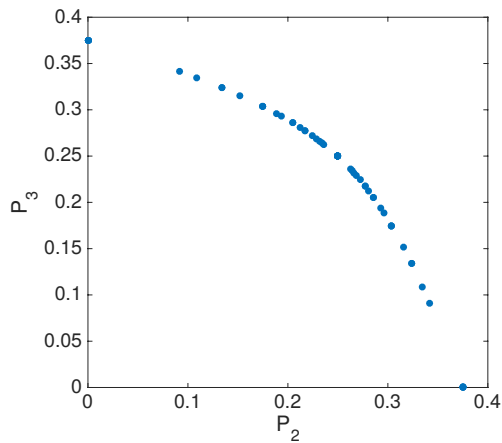}
    \caption{Points on the Pareto-front for the $3$-bus network behind Fig. \ref{fig:jacobian}. They are obtained by solving (\ref{equ:bus_d}) with different  $\bl{z}$'s. The located points are the same as the points on the true boundary in Fig. \ref{fig:jacobian}.}
    \label{fig:paretofront3}
    \end{center}
\end{figure}

Given a network and by varying the search direction $\bl z$, we can use (\ref{equ:bus_d}) to find boundary points. Fig. \ref{fig:paretofront3} shows the result for the network in Fig. \ref{fig:jacobian}. By varying $\bl z$, we successfully locate the boundary points in Fig. \ref{fig:paretofront3}, when compared to exhaustive computations in Fig. \ref{fig:jacobian}.

\section{Conclusion}\label{sec:con}
In general, power flow problems are hard to solve. We propose to use rectangular coordinate, which not only provides an integrated geometric understanding of active and reactive powers simultaneously but also a linear Jacobian matrix for loadability analysis. By using such properties, we proposed three optimization-based approaches for ($1$) loadability verification, ($2$) calculating the margin of an operating point, and ($3$) calculating the boundary points. Numerical results demonstrate the capability of the new method. Future work includes applying the geometric understanding to other analysis in power systems.

\bibliographystyle{IEEEtran}
\bibliography{myBibTex}%

\appendix
\subsection{Derivation of \eqref{equ:pfqf_rec}} \label{app:der}
The active and reactive powers in the polar coordinate are
\begin{align*}
p_d &= \sum_{k=1}^{n}|v_d||v_k|\left(g_{dk}\cos\theta_{dk} + b_{dk}\sin\theta_{dk}\right),\\
q_d &= \sum_{k=1}^{n}|v_d||v_k|\left(g_{dk}\sin\theta_{dk} -b_{dk}\cos\theta_{dk}\right).
\end{align*}

Let $v_{d,r} = |v_d|\cos \theta_d$ and $v_{d,i} = |v_d|\sin \theta_d$, we have
\begin{align}
p_{d} =& v_{d,r}\sum_{k\in  \mathcal{N}(d)}v_{k,r}g_{kd}+v_{d,i}\sum_{k\in  \mathcal{N}(d)}v_{k,i}g_{kd} \nonumber \\
&-v_{d,r}^2\sum_{k\in  \mathcal{N}(d)}g_{kd} - v_{d,i}^2\sum_{k\in  \mathcal{N}(d)}g_{kd} \nonumber\\
&+v_{d,i}\sum_{k\in  \mathcal{N}(d)}v_{k,r}b_{kd} - v_{d,r}\sum_{k\in  \mathcal{N}(d)}v_{k,i}b_{kd} \nonumber\\
=&\Big[-\sum_{k\in  \mathcal{N}(d)}g_{kd}\cdot v_{d,r}^2+\sum_{k\in  \mathcal{N}(d)}v_{k,r}g_{kd}\cdot v_{d,r} \nonumber\\
-& \sum_{k\in  \mathcal{N}(d)}v_{k,i}b_{kd}\cdot v_{d,r}\Big]+\Big[\sum_{k\in  \mathcal{N}(d)}v_{k,i}g_{kd}\cdot v_{d,i} \nonumber\\
+&\sum_{k\in  \mathcal{N}(d)}v_{k,r}b_{kd}\cdot v_{d,i}- \sum_{k\in  \mathcal{N}(d)}g_{kd}\cdot v_{d,i}^2\Big] \nonumber\\
=&t_{d,1}\cdot v_{d,r}^2 + t_{d,2}\cdot v_{d,r}+ t_{d,1}\cdot v_{d,i}^2 + t_{d,3}\cdot v_{d,i}. \label{eqn:pf_2}
\end{align}
Similarly, expanding the equation for $q_d$ and collecting terms gives
\begin{align}
  q_d=&\left(v_{d,r}^2+v_{d,i}^2\right)\sum_{k\in \mathcal{N}(d)}b_{kd} \\&-\left[\sum_{k\in \mathcal{N}(d)}(v_{k,r}b_{kd}+ v_{k,i}g_{kd})\right]  v_{d,r} \nonumber \\
  &+ \left[\sum_{k\in \mathcal{N}(d)}(v_{k,r}g_{kd}- v_{k,i}b_{kd})\right] v_{d,i} \nonumber \\
  =& t_{d,4} v_{d,r}^2-t_{d,3}v_{d,r}+t_{d,4}v_{d,i}^2+t_{d,2}v_{d,i}, \label{eqn:qf_2}
\end{align}
where $t_{d,1},t_{d,2},t_{d,3},t_{d,4}$ are given in \eqref{equ:t_def}.

\subsection{Proof of Lemma \ref{lem:circle}}
\label{lem:circle_app}
From \eqref{eqn:pf_2}, we can complete the square to write it as an equation of a circle
\begin{equation*}
\frac{p_d}{t_{d,1}}+\frac{t_{d,2}^2}{4t_{d,1}^2}+\frac{t_{d,3}^2}{4t_{d,1}^2}=\left(v_{d,r}+\frac{t_{d,2}}{2t_{d,1}}\right)^2+\left(v_{d,i}+\frac{t_{d,3}}{2t_{d,1}}\right)^2.
\end{equation*}
Therefore, the center of the circle is $\left(-\frac{t_{d,2}}{2t_{d,1}},-\frac{t_{d,3}}{2t_{d,1}}\right)$ with the radius being the square root of  $\frac{p_d}{t_{d,1}}+\frac{t_{d,2}^2}{4t_{d,1}^2}+\frac{t_{d,3}^2}{4t_{d,1}^2}$. Interestingly, since $t_{d,1}=-\sum_{k\in  \mathcal{N}(d)}g_{kd}$ is always negative, as the load $p_d$ increases, the radius approaches $0$. This places an upper bound on the possible active power consumption at a bus.

Similarly, we can write \eqref{eqn:qf_2} as
\begin{equation*}
\frac{q_d}{t_{d,4}}+\frac{t_{d,3}^2}{4t_{d,4}^2}+\frac{t_{d,2}^2}{4t_{d,4}^2}=\left(v_{d,r}-\frac{t_{d,3}}{2t_{d,4}}\right)^2+\left(v_{d,i}+\frac{t_{d,2}}{2t_{d,4}}\right)^2,
\end{equation*}
which is a circle with center $\left(\frac{t_{d,3}}{2t_{d,4}},-\frac{t_{d,2}}{2t_{d,4}}\right)$ with the radius being the square root of $\frac{q_d}{t_{d,4}}+\frac{t_{d,3}^2}{4t_{d,4}^2}+\frac{t_{d,2}^2}{4t_{d,4}^2}$. Again, since $t_{d,4}=\sum_{k\in\mathcal{N}(d)} b_{kd}$ is always negative (assuming transmission lines are inductive), the radius decreases as $q_d$ increases. 

\subsection{Incorporating Practical Constraints into Pareto-Front Method}\label{sec:add}
\noindent\textbf{Adding Active Power Constraints.} Box and linear constraints on the active power can be easily added to \eqref{opt:verify} by restricting the direction of movements. Given an operating point, for the constraints on active power that are tight, we can add constraints into \eqref{opt:verify} such that $\bl y$ must move the active powers to stay within the constraint.

\noindent\textbf{Adding Reactive Power Constraints.} The reactive power limits can be visualized in our approach by adding a linear constraint in our algorithms.
	Similar to the definition of $\bd h_d$, let $\bd g_d=[\pd{q_d}{v_{1,r}}, \pd{q_d}{v_{1,i}}, \pd{q_d}{v_{2,r}}, \cdots , \pd{q_d}{v_{n,i}}]^T \in R^{2n}$ be the gradient of $q_d$ with respect to all the state variables, that is, the transpose of the $(n+d)^{th}$ row of the Jacobian matrix.


	To account for the limit of the reactive power constraint, suppose bus $k$ is at its reactive limit, then we modify \eqref{opt:verify} to be
  \begin{subequations} \label{opt:verify_q}
  \begin{align}
    \min_{\bd y} \;\; & 1 \\
    \mbox{s.t. } & \bd y^T \bd h_d \geq 0, \mbox{ for all } d=1,\cdots,n , \\
    & \bd y^T \bd g_k \leq 0, \mbox{ for bus $k$ at reactive power limit}, \label{eqn:limit_q} \\
    & \sum_{d=1}^n \bd y^T \bd h_d =1, 
  \end{align}
  \end{subequations}
where \eqref{eqn:limit_q} specifies that reactive power cannot increase at the buses that hit the limit.

\begin{figure}[ht]
  \centering
  \includegraphics[scale=0.15]{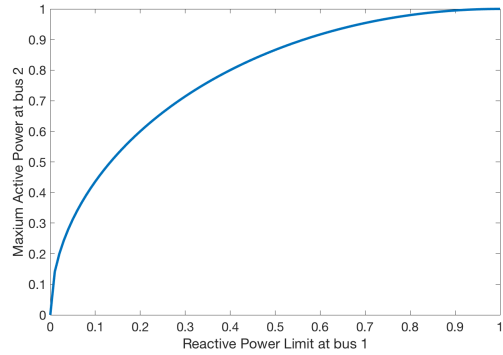}
  \caption{Reactive power limit vs. maximum active power transfer. All values at per unit.}
  \label{fig:twobus_reactive}
\end{figure}

\begin{figure}[ht]
  \centering
  \includegraphics[scale=0.2]{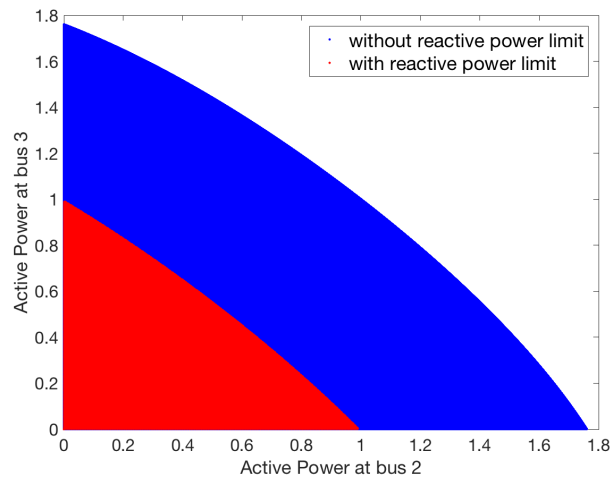}
  \caption{Comparison of feasible power flow region with and without reactive limit on the generator bus for a $3$-bus network.}
  \label{fig:threebus_reactive}
\end{figure}

To visually observe the impact of reactive power limits, we start with a $2$-bus network. We assume that bus $1$ is the slack bus with a reactive limit, and bus $2$ is the PQ load bus. The reactive power limit in this case simply puts a limit on the amount of active power that can be transferred from bus $1$ to bus $2$. Fig.~\ref{fig:twobus_reactive} plots the active power transfer limit at bus $2$ as a function of the reactive limit at bus $1$.
For a $3$-bus system, \eqref{opt:verify_q} also applies. Again, suppose that bus $1$ is the generator (also slack) bus and is limited by reactive power.  For visualization purposes, we plot the feasible region of the achievable powers at bus $2$ and bus $3$ in Fig.~\ref{fig:threebus_reactive}.

\noindent\textbf{Adding Voltage and Current Constraints.} In addition to power constraints, other constraints exist in the system and may become binding before loadability limits are reached. In our approach, incorporating current and voltage limits is straightforward. For example, in the problem in this subsection, we are interested in checking whether a voltage operating point is on the boundary of the feasible region. To include voltage constraints, we can simply add in these as bounds in the voltage space. Similarly, since a current is linear in voltage, current limits can be presented as constraints in the voltage space. After checking these constraints, we can then apply Theorem \ref{thm:boundary} again, which states that checking whether a point is on the boundary can be accomplished by solving a linear program.

\end{document}